\documentclass[12pt]{article}
\usepackage{amssymb}
\usepackage{a4}
\newcommand{\sect}[1]{\setcounter{equation}{0}\section{#1}}

\def\g{\mbox{\bf g\,}}
\def\uqg{\mbox{$U_q{\/\mbox{\bf g}}$ }}
\def\R{\mbox{$\cal R$\,}}
\def\F{\mbox{$\cal F$}}
\def\1{{\bf 1}}
\def\ot{\otimes}
\def\id{\mbox{id}}
\def\A{\mbox{$\cal A$}}
\def\Aq{\mbox{${\cal A}_q$}}

\newcommand{\tr}{\triangleright_q\,}
\newcommand{\trl}{\tilde\triangleright_q\,}
\newcommand{\trc}{\triangleright}
\def\up{\uparrow}
\def\b#1{{\mathbb #1}}
\begin{document}

\title{On $q$-Deformations of Clifford Algebras$^*$}

\author{        Gaetano Fiore \\\\
        \and
        Dip. di Matematica e Applicazioni, Fac.  di Ingegneria\\ 
        Universit\`a di Napoli, V. Claudio 21, 80125 Napoli
        \and
        I.N.F.N., Sezione di Napoli,\\
        Mostra d'Oltremare, Pad. 19, 80125 Napoli
        }
\date{}

\maketitle
\abstract{Several Clifford algebras that are covariant under the 
action of a Lie algebra $\g$ can be deformed in a way consistent
with the deformation of $U\g$ into a quantum group (or into a 
triangular Hopf algebra) $\uqg$, i.e. so as to remain covariant 
under the action of $\uqg$. In this report, after recalling 
these facts, we review our results regarding the formal 
realization of the 
elements of such ``$q$-deformed'' Clifford algebras as
``functions'' (polynomials) in the generators of the 
undeformed ones; in particular, the intruiging interplay 
between the original 
and the $q$-deformed symmetry. Finally, we briefly 
illustrate their dramatic
consequences on the representation theories of the original and 
of the $q$-deformed Clifford algebra, and mention how these results 
could turn out to be useful in quantum physics.}

\vfill
\noindent
- Preprint 99-51 Dip. Matematica e Applicazioni, 
Universit\`a di Napoli

\noindent
$*$ Invited talk given at the ``5th International Conference on 
Clifford algebras
and their Applications in Mathematical Physics, 
Ixtapa-Zihuatanejo (Mexico) June-July 1999.

\newpage

\sect{Introduction}

We first introduce the notions of a $q$-deformed 
Clifford algebra and of a
deforming map in an explicit way on a simple example.

Consider the Clifford algebra $\A$ generated by 
$\1,a^{\up},a^{\downarrow},a^+_{\up},a^+_{\downarrow}$ 
fulfilling
the anticommutation relations
\begin{eqnarray}
&&a^i\, a^j +     a^j\, a^i     = 0 \nonumber\\~
&&a^+_i\, a^+_j +     a^+_j\, a^+_i   = 0 \qquad\qquad i,j=
\up,\downarrow\label{ccr2}\\~
&&a^i\, a^+_j+a^+_j\, a^i      = \delta_j^i\1
\nonumber
\end{eqnarray}
When equipped with the $*$-structure
\begin{equation}
(a^i)^\star=a^+_i                        \label{cstar} 
\end{equation}
this becomes the familiar algebra of creation and annihilation 
operators 
of a fermionic system with two modes (e.g. a spin-up and a 
spin-down 
one-particle state).
$\A$ is then covariant under the action
$\trc$ of $su(2)$, i.e. is a $Usu(2)$-module algebra. 
This  means that (\ref{ccr2}) are left 
invariant by the action of $su(2)$, which is
given by the standard defining action of $su(2)$
on the generators $a^i,a^+_i$ and is
extended to the whole $\A$ according to 
linearity and the Leibniz rule
\begin{equation}
x\trc (\alpha\alpha')=(x\trc \alpha)\alpha'+\alpha(x\trc \alpha').
\label{leibniz}
\end{equation}
We shall not assign the $*$-structure for the moment, so 
\A ~will be covariant under the action of $sl(2,\b{C})$.

$\A$ is the simplest algebra one can $q$-deform: The corresponding
$q$-deformed algebra $\Aq$ ($q\in\b{C}-\{0\}$) is generated by
$\1_q,\tilde A^{\up},\tilde A^{\downarrow},\tilde A^+_{\up},
\tilde A^+_{\downarrow}$ fulfilling the quadratic  
anticommutation relations
\begin{eqnarray}
&& \tilde A^{\up}\tilde A^{\up}   = 0 = \tilde A^{\downarrow}
\tilde A^{\downarrow}\nonumber\\~
&& \tilde A^+_{\up}\tilde A^+_{\up}= 0 =\tilde A^+_{\downarrow}
\tilde A^+_{\downarrow}\nonumber\\~
&& \tilde A^{\up}\tilde A^{\downarrow}+q\tilde A^{\downarrow}
\tilde A^{\up}= 0 \nonumber\\~
&& \tilde A^+_{\up}\tilde A^+_{\downarrow}+q^{-1}
\tilde A^+_{\downarrow}
\tilde A^+_{\up}= 0 \nonumber\\~
&& \tilde A^{\up}\tilde A^+_{\downarrow}+q^{-1}\tilde 
A^+_{\downarrow}
\tilde A^{\up}= 0 \label{dcr2}\\~
&& \tilde A^{\downarrow}\tilde A^+_{\up}+q^{-1}\tilde A^+_{\up}
\tilde A^{\downarrow}= 0 \nonumber\\~
&& \tilde A^{\up}\tilde A^+_{\up}+\tilde A^+_{\up}
\tilde A^{\up} 
= \1_q +(q^{-1}-1)A^+_{\downarrow}\tilde 
A^{\downarrow}\nonumber\\~
&& \tilde A^{\downarrow}\tilde A^+_{\downarrow}+
\tilde A^+_{\downarrow}\tilde A^{\downarrow}= \1_q. 
\nonumber
\end{eqnarray}
This algebra was first introduced in ref. \cite{puwo}.
Clearly $\Aq\stackrel{q\to 1}{\longrightarrow} \A$
if we identify $\tilde A^i\to a^i$, $\tilde A^+_i\to a^+_i$ 
and $\1_q\to \1$ in the limit.
Moreover, (\ref{dcr2}), and hence $\Aq$, are covariant 
under the action $\trl$ of the quantum group $U_qsl(2)$.

One can easily show that $\A,\Aq$ have the same 
Poincar\'e series.
This means that they
have the same dimension (sixteen) as vector spaces, and that 
the set of ordered monomials in the generators
\begin{equation}
{\cal B}_q:=\{\1_q,\tilde A^{\up},\tilde A^{\downarrow},
\tilde A^+_{\up},
\tilde A^+_{\downarrow},\tilde A^{\up}\tilde A^{\downarrow},
\tilde A^+_{\downarrow}\tilde A^{\up},...\}
\end{equation}
is a basis of $\Aq$ and becomes a basis ${\cal B}$ of $\A$ 
in the limit
$q\to 1$. 

What is a $q$-deformed Clifford algebra like (\ref{dcr2}) good 
for? Can it be used to describe the same physics 
as its undeformed
counterpart (\ref{ccr2}), e.g. a second-quantized fermionic 
system, or a different (although similar) one? 
What is the role of the $q$-deformed 
symmetries?
Answering these questions will require 
a comparison not 
only
of the algebraic structures of $\A,\Aq$, but also of
their 
representations. 

Note that the Poincar\'e series requirement
has already an immediate consequence in representation theory: 
it amounts to say that the (left or right) regular 
representation of $\A$ [i.e. the one where the carrier
space is the vector space $U$ associated to the algebra 
$\A$, and the elements of
the algebra act on $U$ by (left or right) multiplication] 
and that of $\Aq$ have the same dimension, since the basis
${\cal B}$ of $U$ and the basis ${\cal B}_q$of $U_q$ 
have the same
number of elements. 

In order to compare the
representation theories it is helpful to ask first a related 
question:
{\bf Can we realize the generators $\tilde A^i,\tilde A^+_i$
as functions (polynomials) in $a^i,a^+_j$? And conversely?}

The answer is yes. For instance, an explicit realization of 
$\tilde A^i,\tilde A^+_i$
 fulfilling (\ref{dcr2}) is
 through
the polynomials
\begin{equation}
\begin{array}{rclcrcl}
A^+_{\downarrow} & := & a^+_{\downarrow}&\qquad
\qquad  A^+_{\uparrow} & := & 
[1+(q^{-1}-1)
n^{\downarrow}]a^+_{\uparrow}
=q^{-n^{\downarrow}}a^+_{\uparrow}  \nonumber \\
A^{\downarrow} & := & a^{\downarrow}&\qquad
\qquad A^{\uparrow} & := &
a^{\uparrow}[1+(q^{-1}-1)
n^{\downarrow}]=
a^{\uparrow}q^{-n^{\downarrow}},
\end{array}
\label{lastf}
\end{equation}
where $n^i:=a^ia^+_i$ (with {\it no} sum over $i$),  
$i=\uparrow, \downarrow$. (The last equalities on the right
are based on the identity $(n^i)^2=n^i$).
For $q\neq 0$ the
transformation (\ref{lastf}) is clearly invertible; 
the inverse
transformation allows an explicit realization of 
$a^i,a^+_i$ fulfilling
(\ref{ccr2}) as polynomials $\tilde a^i, \tilde a^+_i$
in $\tilde A^i,\tilde A^+_i$:
\begin{equation}
\begin{array}{rclcrcl}
\tilde a^+_{\downarrow}& := & \tilde A^+_{\downarrow}&\qquad
\qquad  \tilde a^+_{\uparrow} & := & 
[1+(q-1)N^{\downarrow}]
\tilde A^+_{\uparrow}
=q^{N^{\downarrow}} \tilde A^+_{\uparrow}  \nonumber \\
\tilde a^{\downarrow} & := & \tilde A^{\downarrow}&\qquad
\qquad \tilde a^{\uparrow} & := &
\tilde A^{\uparrow}[1+(q-1)N^{\downarrow}]
=\tilde A^{\uparrow}q^{N^{\downarrow}},
\end{array}
\label{inverse}
\end{equation}
where $N^i:=a^ia^+_i$ (with {\it no} sum over $i$),  
$i=\uparrow, \downarrow$. (The last equalities on the right
are based on the identity 
$(N^{\downarrow})^2=N^{\downarrow}$).
In a more abstract language, through the above 
one can define an algebra isomorphism 
$f:\Aq\to\A[[h]]$ ($\A[[h]]\equiv$the algebra
of formal power series in $h=q-1$ with coefficients in 
$\A$), through
\begin{equation}
\begin{array}{l}
f(\tilde A^i)=A^i \qquad\qquad f(\tilde A^+_i)=A^+_i \cr
f(\alpha\beta)=  f(\alpha)f(\beta) \qquad \forall 
\alpha,\beta \in \A
\end{array}
\label{Arealiz}
\end{equation}

Now we illustrate the usefulness of deforming maps to compare
the representation theories of $\A$, $\Aq$. We ask:  
can $f$ be seen also as an {\it operator map} intertwining 
the
representations of $\A$ and of $\Aq$?
In other
words, given a representation of $\A$ 
(resp. \Aq) on a vector space $V$ (resp. $V_q$), 
does $V$ (resp. $V_q$)
carry also a representation of $\Aq$ (resp. \A)? The 
answer is clearly yes,
since $A^i,A^+_j$ (resp. $\tilde a^i,\tilde a^+_i$), 
being polynomials
in $a^i,a^+_i$ (resp. in $\tilde A^i,\tilde A^+_i$), 
result well-defined
operators on  $V$ (resp. $V_q$)\footnote{It is 
interesting to note that
this is {\it not} true instead for $q$-deformations of 
Weyl algebras 
(the algebras obtained by replacing the anticommutators 
in relations
(\ref{ccr2}) by commutators. The latter are in fact 
infinite-dimensional
as vector spaces, and the corresponding realizations 
$A^i,A^+_i$
are formal power series (instead of polynomials)
in $a^i,a^+_i$, so strictly speaking do not belong to 
$\A$ but just to
a suitable completion of $\A$. Correspondingly, it is 
not guaranteed 
that $A^i,A^+_i$ can be defined as operators on the 
the corresponding
vector spaces. In fact, it was e.g. explicitly shown 
\cite{puwo} that
there are many (inequivalent) irreducible 
$*$-representations of the 
simplest $U_qsu(N)$-covariant deformed Weyl 
$*$-algebra $\Aq$, 
whereas just one of the corresponding undeformed 
partner $\A$. In Ref. 
\cite{fiojmp} we re-read this result by showing 
that the corresponding
objects $\tilde a^i,\tilde a^+_i\in\Aq$ become 
ill-defined operators on all but one 
$*$-representation of \Aq. The same
might in principle occur for Clifford algebras 
with an infinite number of
generators}. Thus, the classification of 
the representations
of $\A$ (resp. \Aq) will determine also 
the classification
of the representations of $\Aq$ (resp. \A). 
From deforming maps
one can extract also more specific informations. 
For instance, if we endow
$\A$ also with the star structure (\ref{cstar}), 
the one that is
compatible with (\ref{ccr2})  and the
action of the compact section $Usu(2)$ of $Usl(2)$,  
the corresponding star structure of $\Aq$ compatible 
with (\ref{dcr2})
and the action of $U_qsu(2)$ exists only for real 
$q$ and reads  
\begin{equation}
(\tilde A^i)^{\star_q}=\tilde A^+_i       \label{qstar} 
\end{equation}
It is easy to see that (\ref{lastf}) [resp. 
(\ref{inverse})] allows also a 
realization of $\star_q$ (resp. $\star$) as $\star$  
(resp. $\star_q$). 
Since there is (up to unitary equivalences)
a unique $*$-representation of the $*$-algebra $\A$ and 
a unique $*$-representation of the $*$-algebra $\Aq$
as one can check by direct inspection, we conclude
that they correspond to each other in the above 
identification.

{\bf How have we found (\ref{lastf})? Does it keep track 
of the $U_qsu(2)$
symmetry?}

In this report we shall present a systematic approach, 
based partly
on the
works \cite{fiojmp,fiormp,fiojpa98} (see also 
Ref.'s \cite{fio})
to answer the latter questions for arbitrary 
`$q$-deformed 
Clifford algebras'.  Incidentally, the approach 
works not only for Clifford, but also for $q$-deformed 
Weyl algebras. 

\sect{General framework}

The general setting is the following. The undeformed 
Clifford 
algebra $\A$ is covariant under some Lie algebra $\g$ 
and the 
deformed one $\Aq$ under the quantum group \cite{dr2} 
(or under
a triangular deformation) $\uqg$.
The undeformed algebra $\A$ is generated by 
$\1, a^i,a^+_j$ fulfilling
\begin{eqnarray}
&&a^i\, a^j +     a^j\, a^i     = 0 \nonumber\\~
&&a^+_i\, a^+_j +     a^+_j\, a^+_i   = 0 \label{ccr}\\~
&&a^i\, a^+_j+a^+_j\, a^i      = \delta_j^i\1
\nonumber
\end{eqnarray}
and transforms under the action $\trc$ of $\g$ 
according to some law
\begin{equation}
    x\trc a^+_i=\rho(x)^j_ia^+_j\qquad\qquad
      x\trc a^i=\rho(Sx)_j^ia^j;   \label{pippo}
\end{equation}
here $x\in\g$, $Sx=-x$ and $\rho$ denotes
some matrix representation of \g. Clearly $a^i$ 
transform under 
the contragradient representation of the $a^+_i$ 
one. The action
$\trc$ is extended to all of $U\g\times \A$ imposing 
linearity,
the Leibniz rule (\ref{leibniz}) and the law
\begin{equation}
(xx')\trc \alpha= x\trc (x'\trc\alpha).
\end{equation}
As a consequence, for $x\in U\g$
\begin{equation}
x\trc(\alpha\alpha')
\equiv\sum_i(x^i_{(1)}\trc\alpha)(x^i_{(2)}\trc\alpha'),
\end{equation} 
where the coproduct  of $U\g$ 
$\Delta(x)=\sum_i x^i_{(1)}\otimes x^i_{(2)}$ is 
defined by
$\Delta(x)=x\otimes \1+ \1\otimes x$ for $x\in\g$ and
is extended to all of $U\g$ as an algebra homomorphism.
All this is possible because the action
of $\g$ is manifestly compatible with the 
anticommutation relations 
(\ref{ccr}). Formulae (\ref{pippo}), where now 
we have extended 
$S$ to the whole $U\g$ as the antipode,  
give also the standard extension of $\trc$ to 
$x\in U\g$.

By definition the corresponding q-deformed algebra 
$\Aq$ is generated by 
$\1_q,\tilde A^+_i,\tilde A^i$
(the generators are enumerated by the same index) 
fulfilling deformed anticommutation relations 
with a quadratic structure as in (\ref{ccr}),
of the form
\begin{eqnarray}
&&P^+_q{}_{ij}^{hk}  \tilde A^+_h 
 \tilde A^+_k  
=0 \nonumber\\~
&&P^+_q{}^{ij}_{hk}  \tilde A^k \tilde A^h = 0 
\label{dcr}\\
&& \tilde A^i\tilde A^+_j+ P_q{}^{ih}_{jk}
\tilde A^+_h \tilde A^k= \delta^i_j\,\1_q.
\nonumber
\end{eqnarray}
$P^+_q, P_q$ are matrices with entries in $\b{C}$. 
$P_q$ is a 
$\uqg$-covariant deformation of
the ordinary 
permutator $P$, defined
by $P^{ij}_{hk}=\delta^i_k\delta^j_h$;
$P^+_q$ is the $\uqg$-covariant 
deformation of the ordinary symmetric projector 
$P^+=\frac {1+P}2$,
and is a projector itself, $(P^+_q)^2=P^+_q$. 
Thus, in the limit
$q\to 1$, the $q$-deformed anticommutation relations
(\ref{dcr}) reduce to (\ref{ccr}). Moreover, 
we additionally require
that $\A,\Aq$ have the same
Poincar\'e series.
\uqg-covariance means that (\ref{dcr}) are 
compatible with
the action $\trl$ of $U_q\g$. The latter is defined on the
generators by the law 
\begin{equation}
x\,\trl\tilde A^+_i = \rho_q{}^j_i(x)\tilde A^+_j
\qquad\qquad
x\,\trl\tilde A^i = \rho_q{}^i_j(S_q x)\tilde A^j,
\label{qtrans}
\end{equation}
(here $x\in U_q\g$, $S_q$ is the antipode of $U_q\g$,
$\rho_q$ the quantum group deformation of $\rho$),
and 
is extended on all of $\Aq$ through a modified 
Leibniz rule,
which we shall give below in (\ref{modalg}). It is 
exactly
the requirement of $\uqg$-covariance that determines 
the
form of $P^+_q, P_q$ in (\ref{dcr}).

Except sometimes the case that $q$ is a root of 
unity,
any representation $\rho$ of $U\g$ admits a
$q$-deformation into a representation $\rho_q$ of \uqg 
(in particular this implies that $\rho,\rho_q$ 
have the same dimension), 
and the corresponding projectors $P,P^+$ admit 
\uqg-covariant deformations
$P^+_q, P_q$. However, it is not guaranteed that the  
corresponding 
relations (\ref{dcr}) yield $\A,\Aq$ with the same
Poincar\'e series.

This is guaranteed for {\it arbitrary} $\rho$ only in a 
less general context,
namely if \uqg is a {\it triangular} 
deformation \cite{dr1} of
the Hopf algebra $U\g$ 
(e.g. a Jordanian \cite{ohn}, or a Reshetikin \cite{Resh} 
deformation); 
then $P_q=\hat R:=PR$, $P^+_q=\frac {1+\hat R}2$,
where $R=(\rho_q\otimes\rho_q) \R$
and  $\R$ is a universal object belonging to 
$\uqg\otimes\uqg$
called the universal $R$-matrix \cite{dr1,dr3}. 
Triangularity means
that $\hat R^2=1$. A special case with a broad 
spectrum of potential 
physical applications is when $\rho$ is a direct sum 
$\rho=\bigoplus\limits_{\alpha=1}^M\rho'$ of $M$ 
copies of a simpler
representation $\rho'$; if the latter describes the 
symmetry of
some dynamical system ${\cal S}'$, $\rho$ should 
describe the
symmetry of the composite dynamical system ${\cal S}$ 
obtained taking $M$ copies of ${\cal S}'$.
The $M$ copies could correspond e.g. to different sites 
in some $d$-dimensional lattice, or (if $M=\infty$) 
to different space(time) 
points, resp. in
condensed matter physics or quantum 
field theory. 

If \uqg is a quantum group in the strict sense, 
i.e. a {\it quasitriangular} but not triangular
deformation of $U\g$ \cite{dr2}, the Poincar\'e 
series condition is fulfilled essentially only if
\begin{enumerate}
\item $\g=sl(N) \cite{puwo}, sp(N=2n)$ \cite{fiojpa98} and 
      $\rho=\rho_N\equiv$ $N$-dimensional 
      defining representation
   of $\g$ (e.g. the {\bf N} of $sl(N)$);
\item $\g=sl(N)$ and 
$\rho=\bigoplus\limits_{\alpha=1}^M\rho_N$ 
\cite{fiojpa98} for some integer $M>1$.
\end{enumerate}

In case 1. $P_q=q^{-1}\hat R_N\equiv q^{-1}PR_N$ and the 
projector
$P^+_q$ is given by
\begin{equation}
P^+_q=\cases{\frac{1+q\hat R_N}{q+q^{-1}} \qquad \qquad 
\mbox{ if \g=}sl(N)\cr
\frac{\hat R_N^2+(q^{-1-N}+q^{-1})\hat R_N+ q^{-2-N}1}
{(q+q^{-1})(q-q^{-1-N})}\qquad \qquad\mbox{if \g=}sp(N),\cr} 
\end{equation}
where $R_N=(\rho_{N,q}\otimes\rho_{N,q}) \R$ and  
$\R\in\uqg\otimes\uqg$ is
the socalled universal $R$-matrix \cite{dr2} of \uqg. 
$R_N$ is denoted as
the $R$-matrix of \uqg in the $q$-deformed defining 
representation 
$\rho_{N,q}$.

In case 2., contrary to the triangular case, it turns out 
that the resulting commutation relations between the 
{\it different} copies
automatically order the $M$ copies in a definite way, a 
phenomenon which
we have called a `braided chain' \cite{fiojpa98}; 
consequently 
the only physical lattice
in which it would be reasonable to arrange the 
copies would 
be 1-dimensional.If we use greek indices 
$\alpha,\beta,...=1,...,M$
to enumerate the copies in the prescribed order, then
up to some free normalization
factors (which we omit
for the sake of simplicity) the deformed 
anticommutation relations 
(\ref{dcr}) take the form, 
\begin{eqnarray}
&&\tilde A^+_{\alpha i} \tilde A^+_{\beta j}+ 
q\hat R_N{}_{ij}^{hk}  
\tilde A^+_{\beta h} \tilde A^+_{\alpha k}  =0 
\nonumber\\~
&&\tilde A^{\alpha j} \tilde A^{\beta i}+
q\hat R_N{}^{ij}_{hk}  \tilde A^{\beta k} 
\tilde A^{\alpha h} = 0 
\label{dcrM}
\end{eqnarray}
and either
\begin{equation}
\tilde A^{\alpha i}\tilde A^+_{\beta j}+ 
q^{-1}\hat R_N{}^{ih}_{jk}
\tilde A^+_{\beta h} \tilde A^{\alpha k}= 
\delta^i_j\delta^{\alpha}_{\beta}\,\1_q
\end{equation}
or
\begin{equation}
\tilde A^{\alpha i}\tilde A^+_{\beta j}+ \hat 
R_M^{-1}{}^{\alpha\gamma}_{\beta\delta}
\hat R_N{}^{ih}_{jk}
\tilde A^+_{\gamma h} \tilde A^{\delta k}= 
\delta^i_j\delta^{\alpha}_{\beta}\,\1_q,
\end{equation}
with $\alpha\le \beta$ and $\hat R_M$ the braid 
matrix of $sl(M)$. 
The latter $\Aq$
in fact is covariant not only under $U_qsl(N)$, 
but also under
$U_q(sl(M)\times sl(N))$.

The explicit form of the braid matrix $\hat R_N$ 
of $sl(N)$ is
\begin{equation}
\hat R_N=q\sum\limits_{i=1}^N{}e^i_i\otimes e^i_i+
\sum\limits_{\stackrel{i,j=1}{i\neq j}}^N{}
e^i_j\otimes e^j_i+
(q-q^{-1})\sum\limits_{\stackrel{i,j=1}{i< j}}^N{}
e^i_i\otimes e^j_j
\end{equation}
where $q\in\b{C}-0$ and $e^i_j$ is the matrix 
with all 
vanisihing elements except a 1 at the $i$-th row 
and $j$-th
column.

Our problem can be now formulated more technically 
as follows: 
How to determine all possible 
deforming maps, i.e.
algebra isomorphisms 
(over $\b{C}[[h]]$)
$f:\Aq\to\A[[h]]$\footnote{We recall that  
for any algebra
$B$ $B[[h]]$ denotes the ring of formal 
power series in \-
$h$ with coefficients in $B$}, (h:=q-1) 
for the class of deformed 
Clifford algebras
defined above?

\sect{Construction procedure}

First note that if $\alpha\in\A[[h]]$ is any
element of the form $\alpha=\1+O(h)$ and $f$ 
is a deforming map,
one can obtain a new one $f_{\alpha}$  by the 
inner automorhism
\begin{equation}
f_{\alpha}(\cdot):=\alpha f(\cdot)\alpha^{-1};
\label{inner}
\end{equation}
actually the vanishing of the first Hochschild 
cohomology group 
\cite{gerst} of $\A$ implies that
{\it all} deforming maps can be obtained from one
in this manner. Therefore our problem is reduced to 
finding a particular one, what we are going to 
describe below.

Second, note that given any deforming map $f$ and
using $\trl$ we can draw the solid lines in the diagram 
\begin{equation}
\begin{array}{ccccc}
 \uqg & \times & $\Aq$ & \, 
 \stackrel{\trl}{-\!\!\!-\!\!\!-\!\!\!- \!\!\!-\!\!\!-
\!\!\!\longrightarrow}\, & $\Aq$ \cr
 &\updownarrow& \id \times  f  & &\updownarrow\, f \cr
 \uqg & \times & \A[[h]] & \, 
 \stackrel{\tr}{-\,-\,-\rightarrow}\, & \A[[h]];\cr
\end{array}
\label{diagram}
\end{equation}
we define $\tr$ as the map making 
the diagram commutative (in other words
$\tr:= f\circ\trl\circ(\id\ot f^{-1})$, which 
will realize $\trl$ on $\A[[h]]$). One can easily 
realize
that $\tr\neq \trc$,
since there is no Hopf 
algebra isomorphism $U_q\g\to U\g[[h]]$ \cite{dr3}.
For each $f_{\alpha}$ in (\ref{inner}) one finds 
correspondingly also a different $\tr$, in other words 
by varying $\alpha$ one obtains all pairs $(f,\tr)$. 

Our construction strategy will proceed in the opposite
direction: we shall first determine {\it one} 
particular $\tr$, 
then the corresponding deforming map(s) $f$. 
Actually to define 
the latter 
it suffices to find generators 
$A^i,A^+_j\in\A[[h]]$ 
fulfilling (\ref{dcr}) and the analog of
(\ref{qtrans}), and apply formula (\ref{Arealiz}). 
We shall 
first show an Ansatz for $A^i,A^+_j$ which allows 
to fulfil at 
once the trasformation law (\ref{qtrans}); 
the Ansatz is
based on the properties of the 
``Drinfel'd twist'' \cite{dr3}.  
Then we shall determine in the simplest 
cases the free parameters 
appearing in the Ansatz in such a way 
that the commutation 
relations (\ref{dcr}) become fulfilled.

Let us summarize the elements of our 
construction procedure and
of the notation we shall adopt:
\begin{enumerate}

\item \g, a semisimple Lie algebra if the 
      deformation \uqg 
      we are interested in is triangular, or $sl(N)$, $sp(N)$ 
      if the deformation \uqg we are interested 
      in  is a 
      quantum group. As known, one can associate 
      to $U\g$ a 
      cocommutative Hopf algebra 
      $H\equiv(U\g,\cdot,\Delta,\varepsilon,S)$;  
$\cdot,\Delta,\varepsilon,S$ denote the product, coproduct, 
      counit, antipode.
 We shall use the Sweedler's 
      notation 
      $\Delta(x)\equiv x_{(1)}\ot x_{(2)}$: the rhs 
      stands for a sum 
      $\sum_i x^i_{(1)}\ot x^i_{(2)}$ of different 
      terms, but the symbol
      $\sum_i$ is dropped. We shall denote by 
      $H_q\equiv(\uqg,\bullet,\Delta_q,
      \varepsilon_q,S_q,\R)$ the deformation
      of $H$ we are interested in, respectively a 
      triangular Hopf algebra 
      \cite{dr1} or a quantum group \cite{dr2}. 
      $\bullet,\Delta_q,\varepsilon_q,S_q$ denote
      the deformed product, coproduct, counit, antipode, 
      $\R$ the
      universal $R$-matrix. We shall use the Sweedler's 
      notation 
      (with barred indices)
      $\Delta_q(x)\equiv x_{(\bar 1)}\ot x_{(\bar 2)}$.
\item An algebra isomorphism 
      $\varphi_q:\uqg\rightarrow U\g[[h]]$ over 
      ${\bf C}[[h]]$
      whose existence is proved respectively in 
      Ref. \cite{dr1,dr3}:
      $\varphi_q(x\bullet y)=\varphi_q(x)
      \cdot\varphi_q(y)$.
\item A corresponding Drinfel'd twist\cite{dr1,dr3}, 
i.e. an element  
      $\F\equiv\F^{(1)}\!\ot\!\F^{(2)}\!=\!\1^{\ot^2}\!\!
      +\!O(h)$ of $U\g[[h]]\otimes U\g[[h]]$ such that 
      \begin{equation}
      (\varepsilon\ot \id)\F=\1=(\id\ot \varepsilon)\F,
      \qquad\: \:\Delta_q(a)=(\varphi_q^{-1}\ot 
      \varphi_q^{-1})\big
      \{\F\Delta[\varphi_q(a)]\F^{-1}\big\};
      \end{equation}
      the last formula means that, up to the 
      isomorphism $\varphi_q$,
      $\Delta_q$ is
      related to $\Delta$ by a similarity transformation.
\item $\gamma':=\F^{(2)}\cdot S\F^{(1)}$ and 
      $\gamma:=S\F^{-1(1)}\cdot \F^{-1(2)}$. 
      Up to the isomorphism 
      $\varphi_q$,
      $S_q$ and its inverse are related to $S$ by 
      similarity 
      transformations
      involving resp. $\gamma$ and $\gamma'$.
\item The particular representation $\rho_q$ of $\uqg$ 
fulfilling the
      criteria listed after (\ref{qtrans}), and its 
      classical limit 
      (\ref{pippo}).
\item The generalized Jordan-Schwinger algebra 
homomorphism 
      $\sigma:U\g[[h]]$ $\rightarrow\A[[h]]$, 
      defined on the 
      generators by 
      \begin{equation}
      \sigma(\1_{U\mbox{\footnotesize \bf g}})=
      \1\qquad\qquad
      \sigma(x):=
      \rho(x)^i_ja^+_ia^j
      \end{equation}
      $x\in\g$, and extended to the whole $U\g[[h]]$ 
      as an algebra
      homomorphism, $\sigma(yz)=\sigma(y)\sigma(z)$
      and $\sigma(y+z)=\sigma(y)+\sigma(z)$. It is 
      immediate to
      verify that this extension is possible 
      because $\sigma([x,y])=[\sigma(x),\sigma(y)]$. 
      In the $su(2)$
      $\sigma$ takes the well-known form
      \begin{equation}
      \sigma(j_+)\!=\!a^+_{\uparrow}a^{\downarrow},
      \qquad\qquad
      \sigma(j_-)\!=\!a^+_{\downarrow}a^{\uparrow},
      \qquad\qquad
      \sigma(j_0)\!=\!\frac 12(a^+_{\uparrow}a^{\uparrow}\!-
\!
      a^+_{\downarrow}a^{\downarrow}).
      \label{homo}
      \end{equation}
\item The deformed Jordan-Schwinger algebra homomorphism 
      $\sigma_q:\uqg\rightarrow\A[[h]]$, defined by
      $\sigma_q:=\sigma\circ\varphi_q$.
\item The $*$-structures $*,*_q,\star,\star_q$ 
in $H,H_q,\A,\Aq$, if 
      $\A,\Aq$ are $*$-algebras transforming respectively 
      under the Hopf
      $*$-algebras $H,H_q$ with the compatibility 
      condition
      \begin{equation}
      (x\,\tr a)^{\star_q}=S_q^{-1}(x^{*_q})\tr a^{\star_q}.
      \label{condstar}
      \end{equation}
\end{enumerate}

As anticipated, our first step is to guess a 
realization $\tr$ 
of $\trl$ on $\A[[h]]$, instead of $\Aq$. This 
requires 
fulfilling 
\begin{equation}
(xy)\tr a=x\tr(y\tr a) \qquad\qquad\qquad
x\tr(ab)=(x_{(\bar 1)}\tr a)(x_{(\bar 2)}\tr b)
\label{modalg}
\end{equation}
for any $x,y\in \uqg$, $a,b\in \Aq$; these are the 
conditions 
i.e. characterizing a module algebra [also in the
undeformed case, see formulae (\ref{leibniz})]  
There is a simple way to find such a realization, 
namely by setting
\begin{equation}
x\tr a := \sigma_q(x_{(\bar 1)}) a 
\sigma_q(S_q x_{(\bar 2)});
\label{defprop}
\end{equation}
it is easy to check that (\ref{modalg}) are indeed 
fulfilled 
using the basic axioms characterizing the coproduct, 
counit, 
antipode in a generic Hopf algebra.
The guess has been suggested by the undeformed case, 
where the
same conditions  and realization are obtained
for $U\g,\A,\trc$ if in the two previous formulae we 
just erase the 
suffix ${}_q$ and replace $\Delta_q(x)\equiv 
x_{(\bar 1)}\ot x_{(\bar 2)}$
with the cocommutative coproduct $\Delta(x)\equiv 
x_{(1)}\ot x_{(2)}$.
 
Our second step is to realize elements 
$A^i,A^+_j\in\A[[h]]$ that
transform under the action $\tr$ defined by 
(\ref{defprop}) as 
$\tilde A^i,\tilde A^+_j$
 do under $\trl$ 
[see (\ref{qtrans}], namely
\begin{equation}
x\,\tr A^+_i = \rho_q{}^j_i(x)A^+_j
\qquad\qquad
x\,\tr A^i = \rho_q{}^i_j(S_q x) A^j,.
\label{qqtrans}
\end{equation}
Note that $a^i,a^+_j$ do {\it not} transform
in this way. In Ref. \cite{fiojmp} we proved 
that the following 
objects do: 
\begin{equation}
\begin{array}{lll}
A_i^+ &:= & u\, \sigma(\F^{(1)})a_i^+
\sigma(S \F^{(2)}\gamma)\, u^{-1} \cr
A^i&:= & v\, \sigma(\gamma'S \F^{-1(2)})a^i
\sigma(\F^{-1(1)}) v^{-1};          \cr                
\end{array}
\label{def3}
\end{equation}
the result holds for any choice of 
$g$-invariant elements $u,v=\1+O(h)$ 
in $\A[[h]]$, in particular for $u=v=\1$.

The third step is to fix $u,v$ in such a way 
that the 
deformed commutation relations (\ref{dcr})  
are fulfilled.
One can easily show that the latter may fix at 
most the product 
$u v^{-1}$.
 In the case that $\uqg$ is 
triangular, we showed in
Ref. \cite{fiojmp} that they require 
$u v^{-1}=1$.
In the case that $\uqg$ is the quantum 
group $U_qsl(N)$ and
$\rho=\rho_N$ we proved 
\cite{fiormp} that the deformed commutation 
relations
 require
\begin{equation}
u v^{-1} = \frac{\Gamma(n+1)}{\Gamma_{q^2}(n+1)}.          
\end{equation}
Here $\Gamma$ is Euler's $\gamma$-function, 
$\Gamma_{q^2}$  its 
$q$-deformation characterized by 
$\Gamma_{q^2}(x+1)=(x)_{q^2}
\Gamma_{q^2}(x)$, $n:= \sum_i a^ia^+_i$.
We stress that the above solutions regard the 
case 
of $\rho$ being  the defining representation 
$\rho_N$ of $sl(N)$.
We have yet no  formula yielding the right 
$u v^{-1}$, if any, necessary to fulfil the 
(\ref{dcr}) in the 
other cases.
However it is important to recall that \cite{fiormp} 
in general the (\ref{dcr})  
translate into  conditions on $u v^{-1}$ 
where  the Drinfel'd twist $\F$ appears only through
the socalled `coassociator'
\begin{equation}
\phi:=[(\Delta\ot \id)(\F^{-1})](\F^{-1}\ot \1)
(\1 \ot \F)
[(\id \ot \Delta)(\F).
\label{defphi}
\end{equation}
$\phi$ is known, unlike $\F$, for which up to now
there is an existence proof.
This makes the above conditions explicit and allows 
to search 
$u v^{-1}$ in the general case, if it exists.
The explicit expression for $\phi$ is 
\begin{equation}
\phi=\lim_{x_0,y_0\rightarrow 0^+}
\left\{x_0^{-\hbar t_{12}}
\vec{P}\exp\left[-\hbar\int\limits^{1-y_0}_{x_0}
dx\left({t_{12}\over x}
+{t_{23}\over x-1}\right)\right] 
y_0^{\hbar t_{23}}\right\}
\label{integral}
\end{equation}
where $t=\Delta({\cal C})-\1\otimes {\cal C}-
{\cal C}\otimes\1$, 
${\cal C}$ denoting the quadratic Casimir of $U\g$, 
$t_{12}=t\otimes \1$, $t_{23}=\1\otimes t$,
and  the symbol $\vec{P}$ means that we must understand a
path-ordered integral in the variable $x$.
Note that $\phi=\1^{\ot^3}+ O(h^2)$. (In the 
triangular case, on the
contrary, $\phi\equiv\1^{\otimes^3}$.)

Finally, the residual freedom in the choice of $u,v$ 
is partially fixed
if $H,H_q,\A,\Aq$ are matched (Hopf) $*$-algebras
and we make the additional requirement that $\star$ 
realizes in $\A[[h]]$
the $\star_q$ of $\Aq$. For instance, if 
$(a^i)^{\star}=a^+_i$ 
and $q\in \b{R}^+$ this means
\begin{equation}
(A^i)^{\star}=A^+_i,
\end{equation}
and is fulfilled if we take $u=v^{-1}$.

In general, one can show that the knowledge of 
$(\rho\otimes\id)\F$ 
is sufficient to determine the 
$A^i,A^+_i$ of
formulae (\ref{def3}) 
completely. In the $\g=sl(2)$ case, with $\rho$ 
being the fundamental 
representation,
$(\rho\otimes\id)\F$ is explicitly 
known \cite{zachos2}.
Taking $u=v^{-1}$ one finally finds \cite{fiojmp} 
the result 
(\ref{lastf}) we had anticipated.

Above we have determined in $\A[[h]]$ one 
particular realization 
$A^i,A^+_j$  and $\tr$
of the generators $\tilde A^i,\tilde A^+_j$ and of 
the quantum group 
action.
Its main feature is that 
the $\g$-invariant ground state $|0\rangle$ as well 
as the first excited
states $a^+_i|0\rangle$ of the classical Fock space 
representation
are also respectively $\uqg$-invariant
ground state $|0_q\rangle$ and first excited
states $A^+_i|0_q\rangle$ of the deformed Fock space 
representation.

According to eq. (\ref{inner}) all the other 
realizations are of the form
\begin{equation}
A^{\alpha\,i}=\alpha\,A^i\alpha^{-1}\qquad\qquad
A^{+}_{\alpha\,i}= \alpha\, A^+_i\, \alpha^{-1},
\label{aa'}
\end{equation}
with $\alpha=\1+O(h)\in\A[[h]]$.
They are manifestly covariant under the realization 
$\trc_{h,\alpha}$ of the $\uqg$-action defined by
\begin{equation}
x\,\trc_{h,\alpha}\,a\: :=\:\alpha\;\sigma_q(x_{(\bar 1)}) 
\, a\,\sigma_q(x_{(\bar 2)})\;\alpha^{-1}.
\end{equation}
For these realizations the deformed ground state 
in the Fock space 
representation
reads $|0_q\rangle=\alpha|0\rangle$; if 
$\alpha|0\rangle\neq |0\rangle$ 
the $\g$-invariant ground state and first excited
states of the classical Fock space representation
do not coincide with their deformed counterparts.

\sect{Ordinary vs. $q$-deformed invariants}
\label{inva}
We have introduced two actions on $\A[[h]]$: 
\begin{equation}
\trc: U\g \times \A[[h]]\rightarrow \A[[h]], \qquad\qquad
\tr: \uqg \times \A[[h]]\rightarrow \A[[h]].
\end{equation}
Their respective invariant subalgebras
$\A^{inv}[[h]],\A^{inv}_q[[h]]$ are defined by
\begin{equation}
\A^{inv}_q [[h]] : =\{I\in \A[[h]]\:\: |\:\: x\tr 
I=\varepsilon_q(x) I\qquad\forall x\in \uqg\}
\label{def7}
\end{equation}
and by the analogous equation where all suffices 
${}_q$ are erased.
What is the relation between them? It is easy to 
prove that \cite{fiormp}
\begin{equation}
\A^{inv}_q[[h]]=\A^{inv}[[h]].
\label{prop}
\end{equation}
In other words invariants under the $\g$-action 
$\trc$ are also
$\uqg$-invariants under $\tr$, and conversely, 
although in general
$\g$-covariant objects (tensors) and $\uqg$-covariant 
ones do not
coincide in general!

Let us introduce 
in the vector space $\A^{inv}[[h]]=\A^{inv}_q[[h]]$ 
bases $I^1,I^2,...$ 
and $I^1_q,I^2_q,...$. It is 
immediate
to realize that we can choose the $I^n$ as 
homogeneous,
normal-ordered polynomials in $a^i, a^+_j$ and $I^n_q$ 
as homogeneous,
normal-ordered polynomials in $A^i,A^+_j$, since
$\trc$ acts linearly without changing the degrees in 
$a^i$ and $a^+_j$,
and $\tr$ acts linearly without changing the degrees 
in $A^i$ and $A^+_j$.
Explicitly,
\begin{equation}
\begin{array}{lll}
I^1=a^+_ia^i                &\qquad \qquad &    
I^1_q=A^+_iA^i  \\
I^2=d^{ijk}a^+_ia^+_ja^+_k  &\qquad \qquad &     
I^2_q=D^{ijk}A^+_iA^+_jA^+_k\\
I^3=d'_{kji}a^ia^ja^k       &\qquad \qquad &     
I^3_q=D'_{kji}A^iA^jA^k  \\
I^4=....                    &\qquad \qquad &    
I^4_q=....   
\end{array}
\end{equation}
where the numerical coefficients  $d,d',...$ 
form $\g$-isotropic tensors and the numerical coefficients
$D,D'$ the corresponding $\uqg$-isotropic tensors.
In the quantum group cases considered in Section 2
it  is possible to show that $I^1_q\neq I^1$:
\begin{equation}
I^1_q=\frac{q^{-2 I^1}-1}{q^{-2}-1}              \label{I1}
\end{equation}
 In general 
$I^n_q\neq I^n$, although $I^n_q=I^n+O(h)$.
The propostion (\ref{prop}) implies in particular
\begin{equation}
I^n_q=g^n(\{I^m\},h)=k^n(\{a^i,a^+_j\},h).
\end{equation}
What do the functions $g^n,k^n$ look like?

In Ref. \cite{fiormp} we have found universal formulae 
yielding
the $k^n$'s. The latter turn out to be polynomials in 
$a^i,a^+_i$
of degree higher than the degree in $A^i,A^+_i$ 
[this can be easily
worked e.g. for the invariant $I_q$ given in (\ref{I1})], 
and the degree
difference grows very fast with the number of these 
generators.
It is remarkable that in these universal formulae 
the twist $\F$
appears only through the coassociator $\phi$; 
therefore all the
$k^n$ can be worked out explicitly.

In the case that the Hopf algebra $H_q$ is not a 
genuine quantum group, 
but 
triangular, the coassociator as well as $u,v$ 
are trivial and one finds
$I^n_q=I^n$.

\sect{Final remarks, outlook and conclusions}

We have shown how one can realize a deformed 
$\uqg$-covariant
Clifford algebra $\Aq$
within the undeformed one $\A[[h]]$. Given a 
representation $(\pi,V)$
of $\A$ on a vector space $V$, does it 
provide also a representation
of $\Aq$? In other words, can one interpret 
the elements of $\Aq$ 
as operators acting on $V$, if the elements 
of $\A$ are? If so,
which specific role play the elements 
$A^i,A^+_i$ of $\A[[h]]$?

Repeating the arguments presented in the 
introduction for the
toy-model, one can conlcude
that the answer to the first question is always 
positive, at least
for finite-dimensional Clifford algebras.
In particular, when $q$ is real and the
real structure (\ref{qstar}) is chosen this 
allows to represent
the $q$-deformed Clifford algebra $\Aq$ on the 
standard Fock space
of the original algebra $\A$; in a 
particle-physics interpretation
no exotic statistics are then
involved, but just the ordinary Fermi-Dirac 
characterizing
fermions. Only, $A^i,A^+_j$ do not 
annihilate/create the
undeformed states.

On the other hand quadratic commutation relations
of the type (\ref{dcr}) mean that $A^+_i,A^i$ 
can be interpreted  
as as creators and
annihilators of some excitations; a glance at 
(\ref{def3}), (\ref{aa'})
shows that these are not the undeformed 
excitations,
but some `collective' ones\footnote{The idea that 
deformed
excitations should consist of a compound of ordinary 
ones is not new,
both for fermions and for bosons:
see for instance Ref.'s \cite{vari}}. The last point is:
what could the latter be good for. As an Hamiltonian 
$H$ of the
system we can choose a simple combination of the 
$\uqg$-invariants
$I^n_q$ of section \ref{inva}; the Hamiltonian is 
$\uqg$-invariant and has
a simple polynomial structure in the composite 
operators
$A^i,A^+_j$. $H$ is also $\g$-invariant, but has a 
higher degree
polynomial structure (or more generally a non-polynomial 
structure if $\A$ has an infinite number of generators) 
in the undeformed generators $a^i,a^+_j$. 
This suggests that the 
use of the $A^i,A^+_j$ instead of the $a^i,a^+_j$
should simplify the resolution of the 
corresponding dynamics.

The results presented in the previous paragraphs
could in principle be applied to models in 
quantum field
theory or condensed matter physics by choosing 
representations
$\rho$ which are the direct sum of many copies 
of the same
fundamental representation $\rho_d$; this is what 
we have 
addressed in Ref. \cite{fiojpa98}. The different 
copies could 
correspond respectively to different 
space(time)-points or
crystal sites.

%%%%%%%%%%% Replace the text between here and the next 
%%%%%%%%%% line begining with
%%%%%%%%%%% by the text of your contribution.

%%%%%%%%%%% Insert your bibliography below %%%%%%%%%%

\end{document}